\documentclass[11pt, reqno]{amsart}
\usepackage{amssymb}
\usepackage{amscd}
\usepackage{verbatim,ifthen}
\usepackage{color}
\usepackage{latexsym}
\usepackage{tikz}
\usepackage{tikz-cd}
\usepackage{mathrsfs}
\usepackage{wrapfig}
\usetikzlibrary{shapes}
\usepackage{color}
\usetikzlibrary{arrows.meta}
\usepackage{bbm}
\usetikzlibrary{matrix}
\usetikzlibrary{calc}
\usetikzlibrary{arrows,intersections}
\usepackage{pgfplots}
\usepackage{multicol}
\usepackage{array}
\newcolumntype{M}[1]{>{\centering\arraybackslash}m{#1}}

\usepackage[colorlinks, linkcolor=black, citecolor=magenta, linktocpage,backref=page]{hyperref}
\renewcommand*{\backref}[1]{}
\renewcommand*{\backrefalt}[4]{[{\tiny%
    \ifcase #1 Not cited.%
          \or Cited on page~#2.%
          \else Cited on pages #2.%
    \fi%
    }]}
\addtocontents{toc}{\setcounter{tocdepth}{1}}

\usepackage{amsmath}

\numberwithin{equation}{section}

\let\oldtocsection=\tocsection
 
\let\oldtocsubsection=\tocsubsection

\renewcommand{\tocsection}[2]{\hspace{0em}\oldtocsection{#1}{#2}}
\renewcommand{\tocsubsection}[2]{\hspace{1em}\oldtocsubsection{#1}{#2}}

\def\XXint#1#2#3{{\setbox0=\hbox{$#1{#2#3}{\int}$ }
\vcenter{\hbox{$#2#3$ }}\kern-.6\wd0}}

\usepackage{eucal}
\usepackage{calc}  
\usepackage{enumitem} 
\usepackage{tensor}
\usepackage{graphicx,wrapfig,lipsum}
\usepackage{etoolbox}
\usepackage{marginnote}
\usepackage{lipsum}
\makeatletter
\patchcmd{\@mn@margintest}{\@tempswafalse}{\@tempswatrue}{}{}
\patchcmd{\@mn@margintest}{\@tempswafalse}{\@tempswatrue}{}{}
\reversemarginpar 
\makeatother
\usepackage{scrextend}

\makeatletter
\DeclareRobustCommand\widecheck[1]{{\mathpalette\@widecheck{#1}}}
\def\@widecheck#1#2{%
    \setbox\z@\hbox{\m@th$#1#2$}%
    \setbox\tw@\hbox{\m@th$#1%
       \widehat{%
          \vrule\@width\z@\@height\ht\z@
          \vrule\@height\z@\@width\wd\z@}$}%
    \dp\tw@-\ht\z@
    \@tempdima\ht\z@ \advance\@tempdima2\ht\tw@ \divide\@tempdima\thr@@
    \setbox\tw@\hbox{%
       \raise\@tempdima\hbox{\scalebox{1}[-1]{\lower\@tempdima\box
\tw@}}}%
    {\ooalign{\box\tw@ \cr \box\z@}}}
\makeatother

\title{Second-Order Estimates for Collapsed Limits of Ricci-flat K\"ahler Metrics}
\author{Kyle Broder}
\address{Mathematical Sciences Institute, Australian National University, Acton, ACT 2601, Australia}
\address{BICMR, Peking University, Beijing, 100871, People's Republic of China}
\email{kyle.broder@anu.edu.au}

\usepackage{soul}
\begin{document}

\maketitle

\begin{abstract}
We show that the singularities of the twisted K\"ahler--Einstein metric arising as the long-time solution of the K\"ahler--Ricci flow or in the collapsed limit of Ricci-flat K\"ahler metrics is intimately related to the holomorphic sectional curvature of reference conical geometry. This provides an alternative proof of the second-order estimate obtained by Gross--Tosatti--Zhang \cite{GTZ2019} with explicit constants appearing in the divisorial pole.
\end{abstract}

\section{Introduction}
Declare a compact K\"ahler manifold $(X, \omega)$ to be \textit{Calabi--Yau} if the canonical bundle $K_X$ is holomorphically torsion, i.e., $K_X^{\otimes \ell} \simeq \mathcal{O}_X$ for some $\ell \in \mathbb{N}$. A \textit{fiber space} is understood to mean a surjective holomorphic map $f : X \longrightarrow Y$ from a K\"ahler manifold $X$ onto a normal irreducible reduced projective variety $Y$, with connected, positive-dimensional fibers. A fiber space is said to be \textit{Calabi--Yau} if the smooth fibers are Calabi--Yau in the above sense.  Such Calabi--Yau fiber spaces arise in the study of the K\"ahler--Ricci flow on compact K\"ahler manifolds with semi-ample canonical bundle (where $f$ is given by the Iitaka map $\Phi_{| K_X^{\otimes \ell} |}$), and collapsed limits of Ricci-flat K\"ahler metrics \cite{SongTian2007, SongTian2012, SongTianZhang2019,Tosatti2009,TosattiAdiabaticLimits, Tosatti2011,Tosatti2012,TosattiWeinkoveYang, TosattiZhang2014, TosattiZhangInfiniteTimeSingularities,TosattiZhang, GTZ2013,GTZ2016,GTZ2019,LiTosatti, Broder}.  These families of K\"ahler metrics are known to converge to a twisted K\"ahler--Einstein metric $\omega_{\text{can}}$ on the base of such fiber spaces, satisfying $$\text{Ric}(\omega_{\text{can}}) \ = \ \lambda \omega_{\text{can}} + \omega_{\text{WP}}.$$ Here, $\omega_{\text{WP}}$ is the Weil--Petersson metric measuring the variation in the complex structure of the smooth fibers of $f$.  \\

To state the main results of the present manuscript,  let $\tilde{f} : M^m \longrightarrow N^n$ be a Calabi--Yau fiber space, where we assume that $K_M$ is either holomorphically torsion or semi-ample.  Denote by $\text{disc}(\tilde{f})$ the discriminant locus of $\tilde{f}$, i.e., the set of point $p \in N$ over which the corresponding fiber $\tilde{f}^{-1}(p)$ is singular.  We denote by $\mathcal{D}$ the divisorial component of $\text{disc}(\tilde{f})$, and let $\pi : (Y, \mathcal{E}) \longrightarrow (N,\mathcal{D})$ be a log resolution such that $(Y, \mathcal{E})$ is log smooth.  By making a base-change along $\pi$, we may construct a Calabi--Yau fiber space $f : X^m \longrightarrow Y^n$ such that $Y$ is a smooth projective variety,  and the discriminant locus of $f$ is $\text{disc}(f) = \mathcal{E}$; in particular, $\text{disc}(f)$ has simple normal crossing support.  Decompose $\mathcal{E}$ into irreducible components $\mathcal{E} = \bigcup_{i=1}^{\mu} \mathcal{E}_i$ and let $s_i$ be local defining sections for $\mathcal{E}_i$.  Endow the associated line bundles $\mathcal{O}_Y(\mathcal{E}_i)$ with smooth Hermitian metrics $| \cdot |_{h_i}^2$. \\

The main conjecture concerning the geometry of these twisted K\"ahler--Einstein metrics is the following \cite{SongTianZhang2019,Tosatti2009,TosattiAdiabaticLimits, Tosatti2011,Tosatti2012, TosattiZhang2014, TosattiZhang,GTZ2013,GTZ2016,GTZ2019, LiTosatti,Broder}:

\subsection*{Conjecture 1.1}\label{MainConjecture}
Let $(\mathcal{Z},d_{\mathcal{Z}})$ denote the metric completion of $(N \backslash \mathcal{D}, \text{dist}_{\omega_{\text{can}}})$. Let $$\Phi : (N \backslash \mathcal{D}, \text{dist}_{\omega_{\text{can}}}) \longrightarrow (\mathcal{Z},d_{\mathcal{Z}})$$ be the corresponding local isometric embedding with singular set $S_{\mathcal{Z}} : = \mathcal{Z} - \Phi(N \backslash \mathcal{D})$. Then \begin{itemize}
\item[(i)] the real Hausdorff codimension of $S_{\mathcal{Z}}$ (inside $\mathcal{Z}$) is at least two;
\item[(ii)] $\left( M, \text{dist}_{\omega_t} \right)$ converges to $(\mathcal{Z},d_{\mathcal{Z}})$ in the Gromov--Hausdorff topology;
\item[(iii)] $\mathcal{Z}$ is homeomorphic to $N$.
\end{itemize}

In particular,  (iii) implies that the (collapsed) Gromov--Hausdorff limit for these metrics carries an algebro-geometric structure.  \\

Tosatti--Zhang \cite{TosattiZhang} showed that if the pullback of the twisted K\"ahler--Einstein metric to the birational model $Y$ had (modulo logarithmic poles) at worst conical singularities, then parts (i) and (ii) of \nameref{MainConjecture} followed. Moreover, if $(N, \mathcal{D})$ is log smooth, then part (iii) also follows.  The required estimate was relaxed by the author \cite{Broder}, showing the following:

\subsection*{Theorem 1.2}
Suppose there is a constant $C>0$ and $d \in \mathbb{N}$ such that for any $\varepsilon >0$,  the twisted K\"ahler--Einstein metric affords the estimate \begin{eqnarray}\label{PSEC}
\pi^{\ast} \omega_{\text{can}} & \leq & \frac{C}{| s_{\mathcal{F}} |^{2\varepsilon}} \left( 1 - \sum_{i=1}^{\mu} \log | s_i |_{h_i}^2 \right)^d \omega_{\text{cone}},
\end{eqnarray}

where $\omega_{\text{cone}}$ is a conical K\"ahler metric and $| s_{\mathcal{F}} |^2 : = \prod_j | s_{\mathcal{F}_j} |_{h_j}^2$ is a shorthand, and $\mathcal{F}_j$ denote the $\pi$--exceptional divisors. \\

We refer to \eqref{PSEC} as the conjectural \textit{partial second-order estimate}. The partial second-order estimate with the divisorial pole large and not explicit was proven by Gross--Tosatti--Zhang \cite{GTZ2019}. There, the Aubin--Yau inequality is used, making use of the curvature lower bound on the reference conical metric obtained by Guenancia--Paun \cite{GuenanciaPaun}.  In general, the curvature of the reference conical metric is not bounded below, but is bounded above (see \cite{JMR2016,LinShen,RubinsteinKE}).  In fact, even in the case of a single smooth divisor, the computation in \cite{JMR2016} shows that the holomorphic sectional curvature of the reference conical metric decays to $-\infty$ near the divisor. \\

The main theorem of the present manuscript is to prove the partial second-order estimate with an explicit expression for the divisorial pole:

\subsection*{Theorem 1.3}\label{MainTheorem}
There are uniform constants $C, \delta_0, \Lambda, \varepsilon_1 >0$ and $d \in \mathbb{N}$ such that, for a conical K\"ahler metric $\omega_{\text{cone}}$, we have  \begin{eqnarray}\label{MAINPSE}
\pi^{\ast} \omega_{\text{can}} & \leq & \frac{C}{| s_{\mathcal{F}} |^{2\delta_0(\Lambda + \varepsilon_1)(n-1)}} \left( 1 - \sum_{i=1}^{\mu} \log | s_i |_{h_i} \right)^d \omega_{\text{cone}}.
\end{eqnarray}

\hfill

In contrast with the partial second-order estimate obtained in \cite{GTZ2019}, the divisorial pole is explicit. This will be made clear in the course of the proof (we avoid detailing them here since it becomes rather cumbersome to read).  An important and interesting consequence of the explicit nature of the divisorial pole is the following corollary:

\subsection*{Corollary 1.4}
Suppose that the holomorphic sectional curvature of the reference conical metric $\omega_{\text{cone}}$ is \textit{almost non-positive} in the sense that for any $\varepsilon_0 >0$,  we have $\text{HSC}(\omega_{\text{cone}}) \leq \varepsilon_0$. Then there is some $\varepsilon>0$ such that $$\pi^{\ast} \omega_{\text{can}} \ \leq \ \frac{C}{| s_{\mathcal{F}} |^{2\varepsilon}} \left( 1 - \sum_{i=1}^{\mu} \log | s_i |_{h_i}^2 \right)^d \omega_{\text{cone}}.$$

\subsection*{Remark 1.5}
The proof of the partial second-order estimate makes use of a maximum principle argument applied to a function $\mathcal{Q}$ which we construct.  As a consequence,  it suffices to assume that the holomorphic sectional curvature is almost non-positive at the point where $\mathcal{Q}$ achieves its maximum. In particular, we suspect that by modifying the function we consider here,  the computation can be localized such that the maximum occurs sufficiently close to the divisor, where the holomorphic sectional curvature of $\omega_{\text{cone}}$ will be negative. \\

Combining these results, we see that 

\subsection*{Corollary 1.5}
Suppose the holomorphic sectional curvature of the reference conical metric is almost non-positive.  Then \begin{itemize}
\item[(i)] the real Hausdorff codimension of $S_{\mathcal{Z}}$ (inside $\mathcal{Z}$) is at least two.
\item[(ii)] $\left( M, \text{dist}_{\omega_t} \right)$ converges to $(\mathcal{Z},d_{\mathcal{Z}})$ in the Gromov--Hausdorff topology.
\end{itemize}

Moreover,  if the base of the Calabi--Yau fiber space is smooth and the divisorial component of the discriminant locus has simple normal crossings, then $\mathcal{Z}$ is homeomorphic to $N$. \\

We note that in \cite{GTZ2019}, Gross--Tosatti--Zhang obtained the partial second-order estimate with the divisorial pole being a large constant $A>0$, which is, unfortunately, not explicit.  The present work owes substantially to \cite{GTZ2019},  from which we build upon here. We also note that there has been some recent progress on understanding the Gromov--Hausdorff limit by Li--Tosatti in \cite{LiTosatti}.

\section{Previously known results}
\subsection*{K\"ahler--Ricci flow}
Hamilton's Ricci flow has proven itself to be a natural candidate for a geometric flow which deforms a fixed K\"ahler metric to a canonical metric. Indeed, the Ricci flow is known to preserve the K\"ahler property of the metric; the Ricci flow starting from a K\"ahler metric $\omega_0$ is therefore referred to as the \textit{(normalized) K\"ahler--Ricci flow}: \begin{eqnarray}\label{NORMKRF}
\frac{\partial \omega_t}{\partial t} &=& - \text{Ric}(\omega_t) - \omega_t, \hspace{1cm} \omega_t \vert_{t=0} = \omega_0.
\end{eqnarray} Cao \cite{Cao} showed that starting with any initial reference K\"ahler metric $\omega_0$ on a Calabi--Yau manifold, the K\"ahler--Ricci flow converges smoothly to the Ricci-flat K\"ahler metric in the polarization $[\omega_0]$.  Similar results hold for the K\"ahler--Ricci flow on canonically polarized manifolds, i.e., compact K\"ahler manifolds with ample canonical bundle. For Fano manifolds which admit a K\"ahler--Einstein metric, or more generally, admit a K\"ahler--Ricci soliton, the K\"ahler--Ricci flow converges exponentially fast to the K\"ahler--Einstein metric or K\"ahler--Ricci soliton. From the work of Tian--Zhang \cite{TianZhang2006}, it is known that the maximal existence time $T$ for the K\"ahler--Ricci flow on a compact K\"ahler manifold $X$ is determined by cohomological data: \begin{eqnarray}\label{SharpLocalExistence}
T \ = \ \sup \{ t \in \mathbb{R} : [\omega_0] + t [K_X ] > 0 \}.
\end{eqnarray} In particular, this gives a sharp local existence theorem for the K\"ahler--Ricci flow, and the K\"ahler--Ricci flow encounters finite-time singularities if and only if the flow intersects the boundary of the K\"ahler cone in finite time.  If the canonical bundle is nef, then it follows from \eqref{SharpLocalExistence} that the K\"ahler--Ricci flow exists for all time.  The resulting solutions, in this case, are referred to as \textit{long-time solutions} of the K\"ahler--Ricci flow.

In \cite{SongTian2007, SongTian2012}, Song--Tian outlined a program for the study of the long-time solutions of the K\"ahler--Ricci flow on compact K\"ahler manifolds with semi-ample canonical bundle. The additional assumption of $K_X$ being semi-ample is fruitful since it endows $X$ with the structure of the total space of a Calabi--Yau fiber space. As before, a \textit{fiber space} is understood to mean a surjective holomorphic map $f : X \longrightarrow Y$ with connected fibers from a compact K\"ahler manifold $X$ onto a normal irreducible, reduced, projective variety $Y$. Such a map is said to be a \textit{Calabi--Yau fiber space} if the smooth fibers $X_y:= f^{-1}(y)$ are Calabi--Yau manifolds in the sense that $K_{X_y}$ is holomorphically torsion.

\subsection*{Sequences of Ricci-flat K\"ahler metrics}
In \cite{TosattiAdiabaticLimits}, Tosatti laid the foundational framework, building on \cite{SongTian2007}, for the study of non-collapsed and collapsed limits of Ricci-flat K\"ahler metrics. Here, let $X$ be a Calabi--Yau manifold and denote by $\mathcal{K}$ the K\"ahler cone of $X$. The K\"ahler cone is an open salient convex cone in the finite-dimensional vector space $H^{1,1}(X, \mathbb{R})$. With respect to the metric topology induced by any norm on $\mathcal{K}$, we let $\overline{\mathcal{K}}$ denote the closure of the K\"ahler cone in $H^{1,1}(X, \mathbb{R})$ and denote by $\partial \mathcal{K}$ its boundary. Fix a non-zero class $\alpha_0$ on the boundary of the K\"ahler cone\footnote{This, of course, implies that the dimension of $H^{1,1}(X, \mathbb{R})$ is at least 2. Otherwise, the only class on the boundary is the zero class, in which case everything is trivial.}, i.e., a nef class, which we may assume to have the form $\alpha_0 = f^{\ast}[\omega_Y]$ for some K\"ahler metric $\omega_Y$ on $Y$. Given a polarization $[\omega_0]$ on $M$, we may consider the path \begin{eqnarray}\label{ChoiceOfPath}
\alpha_t = \alpha_0 + t[\omega_0],
\end{eqnarray} for $0 < t \leq 1$. Yau's solution of the Calabi conjecture furnishes a bijection between the Ricci-flat K\"ahler metrics and the points of $\mathcal{K}$. Hence, given that for each $t > 0$, the class $\alpha_t$ is K\"ahler, Yau's theorem endows $\alpha_t$ with a unique Ricci-flat K\"ahler representative $\omega_t$. To understand the behavior of these metrics, therefore, we can vary the complex structure, keeping the cohomological (or symplectic) data fixed, or keep the complex structure fixed and vary the cohomological data; of course, one can also vary both pieces of data, but we will not discuss that here. The former leads to the study of large complex structure limits, which is important in mirror symmetry; but here we will treat only the latter, given its intimate link with the K\"ahler--Ricci flow and canonical metrics. Therefore, we focus on the problem of understanding the behavior of the metrics $\omega_t$ as the cohomology class $\alpha_t$ degenerates to the boundary of the K\"ahler cone. 

The metrics $\omega_t$ are given by solutions to the Monge--Amp\`ere equation \begin{eqnarray}\label{DEGSEQRKF}
\omega_t^m &=& (f^{\ast} \omega_Y + t \omega_X + \sqrt{-1} \partial \overline{\partial} \varphi_t)^m \ = \ c_t t^{m-n} \omega_X^m, \hspace{1cm} \sup_X \varphi_t =0,
\end{eqnarray}

where the constants $c_t$ are bounded away from $0$ and $+\infty$, and converge as $t\longrightarrow 0$. 

Some remarks concerning the choice of path \eqref{ChoiceOfPath} are in order: The reason for considering such a path was initially provided by the fact that this was the path chosen by Gross--Wilson \cite{GrossWilson} in their study of elliptically fibered K3 surfaces with $\text{I}_1$--singular fibers. But the motivation for sticking with such a path, however, is the formal analog with the K\"ahler--Ricci flow (c.f., \eqref{SharpLocalExistence}).  Much of the behavior of the K\"ahler--Ricci flow can be understood from the study of these sequences of Ricci-flat K\"ahler metrics (and vice versa). Hence, it has become standard practice to treat both contexts simultaneously; and this practice will be maintained here.

\subsection*{Twisted K\"ahler--Einstein metrics}
The first systematic approach to the study of these collapsed limits of Ricci-flat K\"ahler metrics (and the K\"ahler--Ricci flow in this setting) was given by Song--Tian \cite{SongTian2007}. They showed that the metrics $\omega_t$ converge to the pullback of a metric $\omega_{\text{can}}$ on the base of the Calabi--Yau fiber space, which satisfies an elliptic equation away from the discriminant locus: \begin{eqnarray}\label{TwistedKE}
\text{Ric}(\omega_{\text{can}}) &=& \lambda \omega_{\text{can}} + \omega_{\text{WP}},
\end{eqnarray}

where $\omega_{\text{WP}}$ is the Weil--Petersson metric measuring the variation in the complex structure of the smooth fibers. A precise description of $\omega_{\text{WP}}$ was given by Tian \cite{TianModuli}, where he showed that $\omega_{\text{WP}}$ is the curvature form of the Hodge metric on $f_{\ast} \Omega_{M/N}$. Let us note that in \eqref{TwistedKE}, $\lambda=0$ for sequences of Ricci-flat K\"ahler metrics, and $\lambda =-1$ for the long-time solution of the K\"ahler--Ricci flow (see \cite{SongTian2007,SongTian2012,TosattiAdiabaticLimits}). In \cite{TosattiAdiabaticLimits}, it was shown that $\omega_t \to f^{\ast} \omega$ in the $\mathcal{C}_{\text{loc}}^{1,\gamma}(M^{\circ})$ topology of K\"ahler potentials for all $0 < \gamma < 1$. In \cite{TosattiWeinkoveYang} this convergence was improved to local uniform convergence and to $\mathcal{C}_{\text{loc}}^{\alpha}(M^{\circ})$ in \cite{HeinTosatti}. If the generic fiber is a (finite quotient) of a torus, or if the family is isotrivial in the sense that all smooth fibers are biholomorphic, the convergence can be improved to the $\mathcal{C}_{\text{loc}}^{\infty}(M^{\circ})$ topology, see \cite{GTZ2013, HeinTosattiRemarks, TosattiZhangInfiniteTimeSingularities} and \cite{HeinTosatti}, respectively.

If $\kappa = n$, then $K_M$ is nef and big. In this case, it is known \cite{TianZhang2006, Tsuji} that the K\"ahler--Ricci flow converges weakly to the canonical metric on $N$, which is smooth on the regular part of $Y$. Further, in \cite{SongKECurrents}, it is shown that the metric completion of the K\"ahler--Einstein metric on $N^{\circ}$ is homeomorphic to $N$. If $0 < \kappa < n$, then $K_M$ is no longer big, and much less is known; it is this case that we consider here. It is expected that the Ricci curvature of the metrics along the flow remains uniformly bounded on compact subsets of $M^{\circ}$ (see, e.g., \cite{TianZhang}). This is known to be the case when the generic fiber is a (finite quotient of a) torus. Recently, in \cite{FongLee}, the higher-order estimates in \cite{HeinTosatti} were used to obtain the uniform bound on $\text{Ric}(\omega_t)$ when the smooth fibers are biholomorphic. Assuming this uniform bound on the Ricci curvature of the metrics along the flow, one can formulate the \nameref{MainConjecture} for the K\"ahler--Ricci flow.

In \cite{TosattiZhang}, Tosatti--Zhang initiated a program to attack \nameref{MainConjecture} by understanding the nature of the singularities of the twisted K\"ahler--Einstein metric $\omega_{\text{can}}$ near the discriminant locus. They showed that if the canonical metric afforded the following $\mathcal{C}^2$--estimate: \begin{eqnarray}\label{conetype}
\pi^{\ast} \omega_{\text{can}} & \leq & C \left( 1 - \sum_{i=1}^{\mu} \log | s_i |_{h_i}^2 \right)^d \omega_{\text{cone}}, 
\end{eqnarray} 
then parts (i) and (ii) of \nameref{MainConjecture} are true. That is, to prove the conjecture, it suffices to show that (modulo some logarithmic poles), the canonical metric is at worst conical near the discriminant locus. For part (iii), Tosatti--Zhang required the additional assumption that the base of the Calabi--Yau fiber space is smooth and the divisorial component of the discriminant locus has simple normal crossings (in particular, they require the resolution $\pi$ to be the identity).

\section{Proof of the main theorem}
Let us recall the set-up from \cite{GTZ2019}: Let $\tilde{f} : M^m \longrightarrow N^n$ be a Calabi--Yau fiber space with discriminant locus $\text{disc}(\tilde{f})$. Let $\mathcal{D}$ denote the divisorial component of $\text{disc}(\tilde{f})$, and let $\pi : (Y,\mathcal{E}) \longrightarrow (N, \mathcal{D})$ be log resolution such that $Y$ is smooth, $\mathcal{E} = \pi^{-1}(\text{disc}(\tilde{f}))$ has simple normal crossings, and $\pi : Y^{\circ}  = Y \backslash \mathcal{E} \longrightarrow N \backslash \mathcal{D}$ is an isomorphism. Let $\tau : X \to M \times_N Y$ be a birational morphism which induces the Calabi--Yau fiber space $f : X \longrightarrow Y$ with $X$ smooth, $f^{-1}(\mathcal{E})$ a divisor with simple normal crossings, and $\text{disc}(f) = \mathcal{E}$.  The twisted K\"ahler--Einstein metric $\omega_{\text{can}}$ \eqref{TwistedKE} is given by $$\omega_{\text{can}} = \omega_N + \sqrt{-1} \partial \overline{\partial} \varphi$$ for some K\"ahler metric $\omega_N$ (understood in the sense of K\"ahler spaces if $N$ is not smooth, see, e.g., \cite{Moishezon}) and a continuous $\omega_N$--plurisubharmonic function $\varphi$.  Let $\pi : (Y, \mathcal{E}) \longrightarrow (N,\mathcal{D})$ be a log resolution as before. The pulled back metric $\pi^{\ast} \omega_{\text{can}}$ then satisfies \begin{eqnarray*}
(\pi^{\ast} \omega_{\text{can}})^n &=& (\pi^{\ast} \omega_N + \sqrt{-1} \partial \overline{\partial} (\pi^{\ast} \varphi))^n.
\end{eqnarray*}

Recall that $\mathcal{E} = \bigcup_{i=1}^{\mu} \mathcal{E}_i$ is a decomposition of $\mathcal{E}$ into irreducible components. For $\alpha_i \in (0,1] \cap \mathbb{Q}$, we may associate a conical metric $\omega_{\text{cone}}$ with cone angle $2\pi \alpha_i$ along $\mathcal{E}_i$. In more detail, $\omega_{\text{cone}}$ is smooth on $Y \backslash \mathcal{E}$, and in any adapted coordinate system, $\omega_{\text{cone}}$ is uniformly equivalent to the model metric $$\sqrt{-1} \sum_{i=1}^k \frac{dz_i \wedge d\overline{z}_i}{| z_i |^{2(1-\alpha_{j_i})}} + \sqrt{-1} \sum_{i=k+1}^n  dz_i \wedge d\overline{z}_i.$$

\subsection*{Proposition 3.1}\label{GTZHODGE}
(\cite[Theorem 2.3]{GTZ2019}). There is a constant $C>0$,  positive integers $d \in \mathbb{N}$, $0 \leq p \leq \mu$, and rational numbers $\beta_i > 0$ ($1 \leq i \leq p$), $0 < \alpha_i \leq 1$ ($p+1 \leq i \leq \mu$), such that on $Y^{\circ}$, \begin{eqnarray*}
C^{-1} \prod_{j=1}^p | s_j |_{h_j}^{2\beta_j} \omega_{\text{cone}}^n \ \leq \ \pi^{\ast} \omega_{\text{can}}^n \ \leq  \ C \prod_{j=1}^p | s_j |_{h_j}^{2\beta_j} \left( 1 - \sum_{i=1}^{\mu} \log | s_i |_{h_i}^2 \right)^d \omega_{\text{cone}}^n,
\end{eqnarray*}

where $\omega_{\text{cone}}$ is a conical K\"ahler metric with cone angle $2\pi \alpha_i$ along components $\mathcal{E}_i$ with $p+1 \leq i \leq \mu$. \\

Consider the smooth nonnegative function $$H \ : = \ \prod_{j=1}^p | s_j |_{h_j}^{2\beta_j},$$ where the product is over those $\pi$--exceptional components (which we denote by $\mathcal{F}$) on which the volume form $\pi^{\ast}\omega_{\text{can}}^{n}$ degenerates. Away from the exceptional divisor $\mathcal{F}$, \begin{eqnarray*}
\sqrt{-1} \partial \overline{\partial} \log(H) &=& -\sum_{j=1}^p \beta_j \Theta^{(\mathcal{F}_j,h_j)},
\end{eqnarray*} where $\Theta^{(\mathcal{F}_j,h_j)}$ is the curvature form of the Hermitian metric $h_j$ on $\mathcal{F}_j$.  Pulling back the Monge--Amp\`ere equation to $Y$ via $\pi$, we get the degenerate, singular Monge--Amp\`ere equation:

\begin{eqnarray}\label{DEGMA}
\pi^{\ast} \omega_{\text{can}}^n &=& (\pi^{\ast} \omega_N + \sqrt{-1} \partial \overline{\partial}(\pi^{\ast}\varphi))^{n} \ \ = \ \ H \psi \frac{\omega_Y^n}{\prod_{i=1}^{\mu} | s_i |_{h_i}^{2(1-\alpha_i)}},
\end{eqnarray}

where $\omega_N$ is the reference form on $N$ and $\omega_Y$ is a reference form on $Y$. From \nameref{GTZHODGE}, the function $$\psi \ : = \ \frac{\pi^{\ast} \omega_{\text{can}}^{n} \prod_i | s_i |_{h_i}^{2(1-\alpha_i)}}{H \omega_{Y}^{n}},$$ which is smooth on $Y^{\circ} \ = \ Y  \backslash E$, has the following growth control: \begin{eqnarray}\label{GROWTH}
C^{-1} \ \leq \ \psi \ \leq \ C \left( 1 - \sum_{i=1}^{\mu} \log | s_i |_{h_i}^2 \right)^d.
\end{eqnarray}

Away from the divisor $\mathcal{E}$, compute \begin{eqnarray}\label{GRAUERT}
\sqrt{-1} \partial \overline{\partial} \log(1/\psi) &=& \text{Ric}(\pi^{\ast} \omega_{\text{can}}) - \text{Ric}(\omega_{Y}) + \sum_i (1-\alpha_i) \Theta^{(\mathcal{E}_i,h_i)} + \sqrt{-1} \partial \overline{\partial} \log(H) \nonumber \\
& \geq & - \text{Ric}(\omega_{Y}) + \sum_i (1-\alpha_i) \Theta^{\mathcal{E}_i,h_i)} - \sum_j \beta_j \Theta^{(\mathcal{F}_j,h_j)}, 
\end{eqnarray}

where the right-hand side of \eqref{GRAUERT} is smooth on all of $Y$. In light of \eqref{GROWTH}, an old result of Grauert and Remmert \cite{GrauertRemmert} tells us that $\log(1/\psi)$ extends to a quasi-plurisubharmonic function with vanishing Lelong numbers on all of $Y$, satisfying \eqref{GRAUERT} everywhere in the sense of currents. Apply Demailly's regularization technique \cite{Demailly1992}, such that we can approach $\log(1/\psi)$ from above by a decreasing sequence of smooth functions $u_j$ which afford the lower bound \begin{eqnarray*}
\sqrt{-1} \partial \overline{\partial} u_j & \geq & - \text{Ric}(\omega_{Y}) + \sum_i (1-\alpha_i) \Theta^{\mathcal{E}_i,h_i)} - \sum_j \beta_j \Theta^{(\mathcal{F}_j,h_j)}- \frac{1}{j} \omega_{Y}, 
\end{eqnarray*}

on all of $Y$. Similarly, since $\sqrt{-1} \partial \overline{\partial} \log(H) \geq - C \omega_{Y}$ (in the sense of currents) on $Y$, we may regularize $\log(H)$ by smooth functions $v_k$ with $\sqrt{-1} \partial \overline{\partial} v_k \geq -C \omega_{Y}$ and $v_k \leq C$ on $Y$. Out of this, we obtain a regularization of \eqref{DEGMA}: $$\omega_{j,k}^{n} \ = \ c_{j,k} e^{v_k -u_j} \frac{\omega_{Y}^{n}}{\prod_i | s_i |_{h_i}^{2(1-\alpha_i)}},$$ where $$\omega_{j,k} \ = \  \pi^{\ast} \omega_0 + \frac{1}{j} \omega_{Y} + \sqrt{-1} \partial \overline{\partial} \varphi_{j,k}$$ is a cone metric (which exists by the main theorems of \cite{JMR2016, GuenanciaPaun}) with the same cone angles as $\omega_{\text{cone}}$. In particular, the metrics $\omega_{j,k}$ are smooth on $Y^{\circ}$ and there is a uniform (independent of $j$) bound on the $L^p$--norm of $\omega_{j,k}^{n}/\omega_{Y}^{n}$, for some $p>1$. By \cite{EGZKE}, it follows that $$\sup_Y | \varphi_{j,k} | \ \leq \ C,$$ and since $e^{v_k -u_j} \longrightarrow H \psi$, we conclude from the stability theorem in \cite{EGZKE} that as $j,k \longrightarrow \infty$, we have $c_{j,k} \longrightarrow 1$, $\varphi_{j,k} \longrightarrow \pi^{\ast} \varphi$ in $\mathcal{C}^0(Y)$, where $\pi^{\ast} \varphi$ solves \eqref{DEGMA}. Fix a reference metric $\omega_Y$ on $Y$ and consider partial regularizations \begin{eqnarray}\label{REGMETRICS}
\omega_j = \pi^{\ast} \omega_{\text{can}} + \frac{1}{j} \omega_Y + \sqrt{-1} \partial \overline{\partial} \varphi_j
\end{eqnarray} given by solutions of the Monge--Amp\`ere equation \begin{eqnarray}\label{MAGREG}
\omega_j^n &=& c_j H e^{-u_j} \frac{\omega_Y^n}{\prod_{i=1}^{\mu} | s_i |_{h_i}^{2(1-\alpha_i)}}.
\end{eqnarray}

\subsection*{Proof of Theorem 1.1}
For positive constants $\delta,  \delta_0, \Lambda_{\delta, j}, \varepsilon, \varepsilon_0, b, C_1>0$, we will apply the maximum principle to the quantity \begin{eqnarray*}
\mathcal{Q} & : = & \log \text{tr}_{\omega_j}(\omega_{\text{cone}}) + \delta \log(H) + ( (\Lambda_{\delta,j} + \varepsilon_0) b - C_1 \delta) \eta  - (\Lambda_{\delta,j} + \varepsilon_0) \varphi_j  \\
&& \hspace*{6cm} - \delta_0(\Lambda_{\delta,j} + \varepsilon_0) \log | s_{\mathcal{E}} |^2 + \varepsilon \log | s_{\mathcal{E}} |^2.
\end{eqnarray*}

The constants will be described throughout the course of the proof.  We start with the following observation: From \cite[Proposition 5.1]{GTZ2019}, there is a $j$--dependent constant $C_j$ such that $\text{tr}_{\omega_{\text{cone}}}(\omega_j) \leq C_j$.  From the complex Monge--Amp\`ere equation, this implies that for any $\delta \in (0,1)$, \begin{eqnarray*}
\text{tr}_{\omega_j}(\omega_{\text{cone}}) & \leq & \frac{C_j^{n-1}}{c_j H e^{-u_j}} \ \leq \ \frac{C_j^{n-1}}{c_j H^{\delta} e^{-u_j}}.
\end{eqnarray*}

In particular,  there is a $j$--dependent constant (abusively denoted by the same symbol) $C_j$ such that $$H^{\delta} e^{-u_j} \text{tr}_{\omega_j}(\omega_{\text{cone}}) \ \leq \  C_j.$$ This will allow us to apply the maximum principle to terms involving $\text{tr}_{\omega_j}(\omega_{\text{cone}})$ so long as we scale $\text{tr}_{\omega_j}(\omega_{\text{cone}})$ by $H^{\delta} e^{-u_j}$ and add a term which decays to $-\infty$ along the divisor $\mathcal{E}$ (at an arbitrary rate).  This latter term is given by $\varepsilon \log | s_{\mathcal{E}} |^2$, where $0 < \varepsilon \leq \frac{1}{j}$.  To compute $\Delta_{\omega_j} \log \text{tr}_{\omega_j}(\omega_{\text{cone}})$ we observe that the complex Monge--Amp\`ere equation gives the following control of the Ricci curvature of $\omega_j$: \begin{eqnarray*}
\text{Ric}(\omega_j) & = & - \sqrt{-1} \partial \bar{\partial} \log(H)+ \sqrt{-1} \partial \bar{\partial}u_j + \text{Ric}(\omega_Y) + \sqrt{-1} \partial \bar{\partial} \log \prod |s_i |_{h_i}^{2(1-\alpha_i)} \\
& \geq & - \frac{1}{j} \omega_Y \ \geq \ - \frac{C}{j} \omega_{\text{cone}}.
\end{eqnarray*}

Let $g,h$ be metrics underlying $\omega_j$ and $\omega_{\text{cone}}$, respectively.  Then for any holomorphic map $f : (Y^{\circ}, \omega_j) \to (Y^{\circ}, \omega_{\text{cone}})$,  with $(f_i^{\alpha})$ denoting the local expression for $\partial f$, we have (see, e.g., \cite{BroderSBC, BroderSBC2, YangZheng}) \begin{eqnarray*}
\Delta_{\omega_j} | \partial f |^2 & = & | \nabla \partial f |^2 + \text{Ric}_{k \overline{\ell}}^g g^{k \overline{q}} g^{p \overline{\ell}} h_{\alpha \overline{\beta}} f_p^{\alpha} \overline{f_q^{\beta}} - R_{\alpha \overline{\beta} \gamma \overline{\delta}} g^{i \overline{j}} f_i^{\alpha} \overline{f_j^{\beta}} g^{p \overline{q}} f_p^{\gamma} \overline{f_q^{\delta}}.
\end{eqnarray*}

The lower bound on the Ricci curvature of $g$ implies that \begin{eqnarray*}
\text{Ric}_{k \overline{\ell}}^g g^{k \overline{q}} g^{p \overline{\ell}} h_{\alpha \overline{\beta}} f_p^{\alpha} \overline{f_q^{\beta}} & \geq & - \frac{C}{j} h_{\gamma \overline{\delta}}  g^{k \overline{q}} g^{p \overline{\ell}} h_{\alpha \overline{\beta}} f_p^{\alpha} \overline{f_q^{\beta}}.
\end{eqnarray*}

Choose the frame such that at the point where we are computing, we have $f_i^{\alpha} = \lambda_i \delta_i^{\alpha}$, $g_{i\overline{j}} = \delta_{ij}$ and $h_{\alpha \overline{\beta}} = \delta_{\alpha \overline{\beta}}$.  Then \begin{eqnarray*}
- \frac{C}{j} h_{\gamma \overline{\delta}} f_k^{\gamma} \overline{f_{\ell}^{\delta}}  g^{k \overline{q}} g^{p \overline{\ell}} h_{\alpha \overline{\beta}} f_p^{\alpha} \overline{f_q^{\beta}} &=& - \frac{C}{j} \delta_{\gamma \delta} \lambda_k \delta_k^{\gamma} \lambda_{\ell} \delta_{\ell}^{\delta} \delta^{kq} \delta^{p\ell} \delta_{\alpha \beta} \lambda_p \delta_p^{\alpha} \lambda_q \delta_q^{\beta} \ = \ - \frac{C}{j} \lambda_{\alpha}^4.
\end{eqnarray*}

Since $\sum_{\alpha} \lambda_{\alpha}^4 \leq \left( \sum_{\alpha} \lambda_{\alpha}^2 \right)^2 = | \partial f |^4$, we deduce that \begin{eqnarray*}
\text{Ric}_{k \overline{\ell}}^g g^{k \overline{q}} g^{p \overline{\ell}} h_{\alpha \overline{\beta}} f_p^{\alpha} \overline{f_q^{\beta}} & \geq & - \frac{C}{j} \text{tr}_{\omega_j}(\omega_{\text{cone}})^2
\end{eqnarray*}

From \cite{JMR2016, LinShen}, we know that there is a uniform constant $K>0$ such that $\text{HBC}(\omega_{\text{cone}}) \leq K$. This, in particular, implies that the holomorphic sectional curvature of the cone metric $\omega_{\text{cone}}$ affords the (positive) upper bound $\text{HSC}(\omega_{\text{cone}}) \leq K$ for $K>0$.  Royden's polarization argument \cite{Royden,Lu,BroderSBC2,BroderThesis} implies that \begin{eqnarray*}
R_{\alpha \overline{\beta} \gamma \overline{\delta}} g^{i \overline{j}} f_i^{\alpha} \overline{f_j^{\beta}} g^{p \overline{q}} f_p^{\gamma} \overline{f_q^{\delta}} & \leq & K \text{tr}_{\omega_j}(\omega_{\text{cone}})^2.
\end{eqnarray*}

Hence, by Royden's Schwarz lemma \cite{Royden,Lu,BroderThesis}, we see that \begin{eqnarray*}
\Delta_{\omega_j} \log \text{tr}_{\omega_j}(\omega_{\text{cone}}) & \geq & - \left( \frac{C}{j} + K \right) \text{tr}_{\omega_j}(\omega_{\text{cone}}),
\end{eqnarray*}

at any point away from the divisor $\mathcal{E}$.  Therefore,  \begin{eqnarray*}
\Delta_{\omega_j} (\log \text{tr}_{\omega_j}(\omega_{\text{cone}}) + \delta \log(H) - u_j ) & \geq & - \left( \frac{C}{j} + K \right) \text{tr}_{\omega_j}(\omega_{\text{cone}}) \\
&& \hspace*{2cm} + \delta \text{tr}_{\omega_j}(\sqrt{-1} \partial \bar{\partial} \log(H)) - \Delta_{\omega_j} u_j,
\end{eqnarray*}

at any point away from the divisor $\mathcal{E}$.  Since $\sqrt{-1} \partial \bar{\partial} \log(H) \geq - C_1 \omega_Y$,  and adding $u_j$ to the quantity being differentiated (which will not affect the maximum principle), we see that \begin{eqnarray*}
\Delta_{\omega_j} (\log \text{tr}_{\omega_j}(\omega_{\text{cone}}) + \delta \log(H)) & \geq & - \left( \frac{C}{j} + K \right) \text{tr}_{\omega_j}(\omega_{\text{cone}})  - \delta C_1  \text{tr}_{\omega_j}(\omega_Y),
\end{eqnarray*}

at any point away from $\mathcal{E}$.  We know that the conical metric is given by $\omega_{\text{cone}} = \omega_Y + \sqrt{-1} \partial \bar{\partial} \eta$. Hence,   \begin{eqnarray*}
\Delta_{\omega_j}(\log \text{tr}_{\omega_j}(\omega_{\text{cone}}) +\delta \log(H) - \delta C_1 \eta + \varepsilon \log | s_{\mathcal{E}} |^2) & \geq & - \left( \frac{C}{j} + K + C_1 \delta \right) \text{tr}_{\omega_j}(\omega_{\text{cone}}),
\end{eqnarray*}

and we adopt the notation $\Lambda_{\delta,j} : = \frac{C}{j} + K + C_1 \delta$.  To ensure the positivity of the coefficient of $\text{tr}_{\omega_j}(\omega_{\text{cone}})$, let $\varepsilon_1>0$ be a positive constant to be determined later. Then \begin{eqnarray*}
\Delta_{\omega_j}( \log \text{tr}_{\omega_j}(\omega_{\text{cone}}) + \delta \log(H) - C_1 \delta  \eta  + \varepsilon \log | s_{\mathcal{E}} |^2 - (\Lambda_{\delta,j} + \varepsilon_1)  \varphi_j) \hspace*{2cm} \\
\geq - \Lambda_{\delta,j} \text{tr}_{\omega_j}(\omega_{\text{cone}}) - (\Lambda_{\delta,j} + \varepsilon_1)n + (\Lambda_{\delta,j} +\varepsilon_1) \text{tr}_{\omega_j}(\pi^{\ast} \omega).
\end{eqnarray*}

The pullback $\pi^{\ast} \omega$ is degenerate in the directions tangent to the divisor.  Let $\delta_0>0$ be the suitably small constant such that $\pi^{\ast} \omega - \delta_0 \Theta^{(\mathcal{F},h)}$ is positive-definite.  Further, we let $b>0$ be the suitably small constant such that $\pi^{\ast} \omega - \delta_0 \Theta^{(\mathcal{F},h)} + b \sqrt{-1} \partial \bar{\partial} \eta$ has the same cone angles as $\omega_{\text{cone}}$. Hence, we can find a constant $c_0>0$ such that $$\pi^{\ast} \omega - \delta_0 \Theta^{(\mathcal{F},h)} + b \sqrt{-1} \partial \bar{\partial} \eta \ \geq \ c_0 \omega_{\text{cone}}.$$ Hence, with $\mathcal{Q}$ defined above, we see that \begin{eqnarray*}
\Delta_{\omega_j} \mathcal{Q} & \geq & - \Lambda_{\delta, j} \text{tr}_{\omega_j}(\omega_{\text{cone}}) - (\Lambda_{\delta,j} + \varepsilon_1) n + (\Lambda_{\delta,j} + \varepsilon_1) \text{tr}_{\omega_j}(\pi^{\ast} \omega - \delta_0 \Theta^{(\mathcal{F},h)} + b \sqrt{-1} \partial \bar{\partial} \eta) \\
& \geq & (c_0(\Lambda_{\delta,j} + \varepsilon_1) - \Lambda_{\delta,j})  \text{tr}_{\omega_j}(\omega_{\text{cone}}) - (\Lambda_{\delta,j} + \varepsilon_1) n.
\end{eqnarray*}

If $c_0 \geq 1$, then $c_0 (\Lambda_{\delta,j} + \varepsilon_1) - \Lambda_{\delta,j} \geq c_0 \varepsilon_0$, and we take $\varepsilon_1>0$ arbitrary. For instance, if we take $\varepsilon_1 = c_0^{-1}$, then \begin{eqnarray*}
\Delta_{\omega_j} \mathcal{Q} & \geq & \text{tr}_{\omega_j}(\omega_{\text{cone}}) - (\Lambda_{\delta,j} + c_0^{-1}) n.
\end{eqnarray*}

By the maximum principle,  at the point where $\mathcal{Q}$ achieves its maximum, we have $$\text{tr}_{\omega_j}(\omega_{\text{cone}}) \ \leq \ n(\Lambda_{\delta,j} + c_0^{-1}) n \ \leq \ C,$$ for some uniform constant $C>0$.  On the other hand, if $0 < c_0 < 1$, then we choose $\varepsilon_1 = c_0^{-1}(1 + \Lambda_{\delta,j}(1-c_0))$. In this case, we see that  \begin{eqnarray*}
\Delta_{\omega_j}(\mathcal{Q}) & \geq & \text{tr}_{\omega_j}(\omega_{\text{cone}}) - c_0^{-1}n(1+\Lambda_{\delta,j}),
\end{eqnarray*}

and by the maximum principle, at the point where $\mathcal{Q}$ achieves its maximum,  $$\text{tr}_{\omega_j}(\omega_{\text{cone}}) \ \leq \ c_0^{-1} n(1+\Lambda_{\delta,j}).$$ Both estimates imply that $\mathcal{Q} \leq C$ for some uniform constant $C>0$.  Hence, we have the estimate \begin{eqnarray*}
\text{tr}_{\omega_j}(\omega_{\text{cone}}) & \leq & \frac{C}{H^{\delta} | s_{\mathcal{F}} |^{2\delta_0(\Lambda_{\delta,j}+\varepsilon_1)}}
\end{eqnarray*}

everywhere on $Y^{\circ}$.  From the complex Monge--Amp\`ere equation, we see that \begin{eqnarray*}
\text{tr}_{\omega_{\text{cone}}}(\omega_j) & \leq & \text{tr}_{\omega_j}(\omega_{\text{cone}})^{n-1} \frac{\omega_j^n}{\omega_{\text{cone}}^n} \ \leq \ \frac{CH e^{u_j}}{| s_{\mathcal{F}} |^{2\delta_0(\Lambda_{\delta,j} + \varepsilon_1)(n-1)} H^{\delta(n-1)}}.
\end{eqnarray*}

Choosing $\delta < \frac{1}{n-1}$, we see that \begin{eqnarray*}
\text{tr}_{\omega_{\text{cone}}}(\omega_j) & \leq & \frac{Ce^{-u_j}}{| s_{\mathcal{F}} |^{2\delta_0(\Lambda_{\delta,j} + \varepsilon_1)(n-1)}},
\end{eqnarray*}

and letting $j \to \infty$, we obtain the partial second-order estimate \begin{eqnarray*}
\pi^{\ast} \omega_{\text{can}} & \leq & \frac{C}{| s_{\mathcal{F}} |^{2\delta_0(\Lambda_{\delta,j} + \varepsilon_1)(n-1)}} \left( 1 - \sum_{i=1}^{\mu} \log | s_i |_{h_i}^2 \right)^d \omega_{\text{cone}}.
\end{eqnarray*}

\subsection*{Proof of Corollary 1.2}
Suppose the holomorphic sectional curvature of $\omega_{\text{cone}}$ is almost non-positive, i.e.,  for any $\varepsilon_0 >0$,  we have $\text{HSC}(\omega_{\text{cone}}) \leq \varepsilon_0$, then \begin{eqnarray*}
2\delta_0(\Lambda_{\delta,j} + \varepsilon_1)(n-1) & \leq & 2\delta_0 \left( \frac{C}{j} + \varepsilon_0 + C_1 \delta \right)(n-1).
\end{eqnarray*}

The constant $\delta \in (0, \frac{1}{n-1})$ can be made arbitrarily small, and therefore, as $j \to \infty$, we have \begin{eqnarray*}
2\delta_0 \left( \frac{C}{j} + \varepsilon_0 + C_1 \delta \right) (n-1) & \to & 2 \delta_0(n-1)(\varepsilon_0 + C_1 \delta).
\end{eqnarray*}

In particular,  if the holomorphic sectional curvature of the conical K\"ahler metric is bounded above by an arbitrarily small positive constant,  then the multiplicity of the divisorial pole can be made arbitrarily small.

\subsection*{Remark 3.2}
The assumption that the holomorphic sectional curvature of the conical metric is arbitrarily small is not as restrictive as one may expect.  In fact,  we suspect that by appropriately localizing the function $\mathcal{Q}$ used in the maximum principle, one should be able to ensure that the holomorphic sectional curvature of $\omega_{\text{cone}}$ is almost non-positive at the point where $\mathcal{Q}$ achieves its maximum. To see this,  consider the case where the divisor $\mathcal{E}$ consists of a single component, locally described by $\{ z_1 =0 \}$. Then the computation in \cite[Appendix]{JMR2016} shows that $$R_{1 \overline{1} 1 \overline{1}} \sim -\alpha^2(\alpha-1)^2 | z |^{-2(2-\alpha)}.$$  In other words, if the function $\mathcal{Q}$ can be tweaked such that the maximum occurs sufficiently close to $\mathcal{E}$, but not on $\mathcal{E}$,  then the holomorphic sectional curvature of $\omega_{\text{cone}}$ should be negative at this point, and the argument goes through.

\subsection*{The Holomorphic Sectional Curvature}
The main theorem, although primarily concerned with \nameref{MainConjecture}, fits into a more general program the author has initiated.  There are the three primary aspects of the holomorphic sectional curvature: \begin{itemize}
\item[(i)] Symmetries -- The holomorphic sectional curvature of a K\"ahler metric (or more generally, K\"ahler-like metric) determines the curvature tensor entirely.
\item[(ii)] Value distribution of curves -- A compact K\"ahler manifold with negative holomorphic sectional curvature is Kobayashi hyperbolic (every entire curve $\mathbb{C} \to X$ is constant). On the other hand, a compact K\"ahler manifold with positive holomorphic sectional curvature is rationally connected (any two points lie in the image of a rational curve $\mathbb{P}^1 \to X$) \cite{YangHSCYau}.
\item[(iii)] Singularities -- The holomorphic sectional curvature influences the singularities of the geometry.
\end{itemize}

The first facet of the holomorphic sectional curvature is well-known to all complex geometers.  The assertions in statement (ii) are also well-known, but it is worth emphasizing a particular subtlety concerning (ii): Let $(X, \omega)$ be a Hermitian manifold with a Hermitian metric of negative (Chern) holomorphic sectional curvature ${}^c \text{HSC}(\omega) \leq - \kappa <0$. Then $X$ is Brody hyperbolic in the sense that every entire curve $\mathbb{C} \to X$ is constant. The proof is an elementary application of the Schwarz lemma. Suppose there is a non-constant map $f : \mathbb{C} \to X$. Then $\Delta_{\omega_{\mathbb{C}}} | \partial f |^2 \geq - {}^c R(\partial f, \bar{\partial f}, \partial f, \bar{\partial f})$ and by the maximum principle $f$ is constant.  The key point here, however, is that for holomorphic curves, there is only one direction to consider. For holomorphic maps of higher rank (e.g., $f : \mathbb{C}^2 \to X$) it is not clear whether the holomorphic sectional curvature of an arbitrary Hermitian metric gives any control.  This is the reason for the introduction of the real bisectional curvature \cite{YangZheng,LeeStreets} and the Schwarz bisectional curvatures \cite{BroderSBC,BroderSBC2}.

At this point, the reader may be perplexed,  since we used the Schwarz lemma in the proof of the main theorem and this involves a map of rank $>1$. But here,  we note that Royden's polarization argument \cite{Royden,BroderSBC2,BroderThesis} is used.  In particular,  for holomorphic maps into K\"ahler manifolds, the holomorphic sectional curvature gives suitable control.  In other words, we see the Schwarz lemma as an incarnation of aspects (i) and (ii).

The third aspect is far less understood (emphasized by the vague nature of the statement in (iii)). The main theorem of the present manuscript provides evidence for the relationship,  but it is certainly not the only evidence.  Recall that Demailly \cite{Demailly} constructed an example of a projective Kobayashi hyperbolic surface that does not admit a Hermitian metric of negative Chern holomorphic sectional curvature. The construction relies on the following algebraic hyperbolicity criterion \cite{Demailly}:

\subsection*{Theorem}
(Demailly).  Let $(X, \omega)$ be a compact Hermitian manifold with ${}^c \text{HSC}_{\omega} \leq \kappa_0$ for some $\kappa_0 \in \mathbb{R}$.  If $f : \mathcal{C} \to X$ is a non-constant holomorphic map from a compact Riemann surface $\mathcal{C}$ of genus $g$, then \begin{eqnarray*}
2g-2 & \geq & -\frac{\kappa_0}{2\pi} \deg_{\omega}(\mathcal{C}) - \sum_{p \in \mathcal{C}} (m_p-1).
\end{eqnarray*}

\hfill

We note that in \cite{Demailly} (see also \cite{Diverio}) the constant is assumed to be non-positive, but the formula holds more generally (see, e.g., \cite{BroderThesis}).  The example is constructed as a fibration (with hyperbolic base and fiber) with a fiber sufficiently singular to violate the above algebraic hyperbolicity criterion.

\end{document}